\RequirePackage[l2tabu,orthodox]{nag}
\documentclass
[11pt,letterpaper]
{article} 

\usepackage[notes=true,later=false,camera=false]{dtrt}
\usepackage[utf8]{inputenc}
\usepackage{etex}
\usepackage{thmtools}
\usepackage{ stmaryrd }
\usepackage{xspace,enumerate}
\usepackage[T1]{fontenc}
\usepackage[full]{textcomp}
\usepackage[american]{babel}
\usepackage{mathtools}
\usepackage{amsthm}

\usepackage{empheq}

\hypersetup{
colorlinks=true,
urlcolor=Cerulean,
linkcolor=RoyalBlue,
citecolor=OliveGreen,
linktocpage=true,
}
\renewcommand*\backref[1]{\ifx#1\relax \else (pg. #1) \fi}

\usepackage[capitalise,nameinlink]{cleveref}
\crefname{lemma}{Lemma}{Lemmas}
\crefname{fact}{Fact}{Facts}

\crefname{theorem}{Theorem}{Theorems}
\crefname{corollary}{Corollary}{Corollaries}
\crefname{claim}{Claim}{Claims}
\crefname{example}{Example}{Examples}
\crefname{algorithm}{Algorithm}{Algorithms}
\crefname{problem}{Problem}{Problems}
\crefname{definition}{Definition}{Definitions}
\crefname{equation}{Eq.}{Eq.}
\crefname{strategy}{Strategy}{Strategies}

\usepackage{paralist}
\usepackage{turnstile}
\usepackage[framemethod=TikZ]{mdframed}
\mdfsetup{frametitlealignment=\center}
\usepackage{tikz}
\usepackage{caption}
\DeclareCaptionType{Algorithm}
\usepackage{newfloat}
\usepackage{aliascnt}
\newtheorem{theorem}{Theorem}[section]
\newaliascnt{lemma}{theorem}
\newtheorem{lemma}[lemma]{Lemma}
\aliascntresetthe{lemma}\newtheorem*{lemma*}{Lemma}
\newaliascnt{remark}{theorem}
\newtheorem{remark}[remark]{Remark}
\aliascntresetthe{remark}\newtheorem*{remark*}{Remark}
\newaliascnt{definition}{theorem}
\newtheorem{definition}[definition]{Definition}
\aliascntresetthe{definition}\newtheorem*{definition*}{Definition}

\newaliascnt{proposition}{theorem}
\newtheorem{proposition}[proposition]{Proposition}
\aliascntresetthe{proposition}
\newaliascnt{fact}{theorem}
\newtheorem{fact}[fact]{Fact}
\aliascntresetthe{fact}
\newaliascnt{corollary}{theorem}

\aliascntresetthe{corollary}

\newtheorem{algorithm-thm}[theorem]{Algorithm}

\usepackage[
letterpaper,
top=1.2in,
bottom=1.2in,
left=1in,
right=1in]{geometry}
\usepackage{newpxtext} 
\usepackage{textcomp} 
\usepackage{mathpazo}
\usepackage[scr=rsfso]{mathalfa}
\usepackage{bm} 
\let\mathbb\varmathbb
\usepackage{microtype}

\usepackage{footnotebackref}

\allowdisplaybreaks
\newcommand{\FormatAuthor}[3]{
\begin{tabular}{c}
#1 \\ {\small\texttt{#2}} \\ {\small #3}
\end{tabular}
}

\newcommand{\R}{{\mathbb R}}
\newcommand{\N}{{\mathbb N}}

\newcommand{\eps}{\varepsilon}

\newcommand{\F}{{\mathbb F}}

\newcommand{\E}{{\mathbb E}}
\newcommand{\cE}{{\mathcal E}}

\newcommand{\1}{\mathbf{1}}
\newcommand{\ip}[1]{\langle #1 \rangle}

\newcommand{\C}{\mathbb C}

\newcommand{\cH}{\mathcal H}

\newcommand{\cA}{\mathcal A}
\newcommand{\cB}{\mathcal B}
\newcommand{\cP}{\mathcal P}

\newcommand{\cS}{\mathcal S}

\newcommand{\mper}{\,.}
\newcommand{\mcom}{\,,}

\newcommand{\cC}{\mathcal C}

\renewcommand{\emptyset}{\varnothing}
\renewcommand{\geq}{\geqslant}
\renewcommand{\leq}{\leqslant}
\renewcommand{\epsilon}{\varepsilon}

\usepackage{macros}
\title{A Short Proof and Generalization of the Rainbow Cycle Theorem and Improved Bounds for Binary Even-Query LCCs}
\author{\FormatAuthor{Andrew D.\ Lin}{andrewlin@princeton.edu}{Princeton University}
  }
  \date{\today}
\begin{document}

\newcommand{\fhj}{F_\cH(j)}
\newcommand{\fgj}{F_{\gnkt}(j)}
\newcommand{\lcs}{\text{LCS}}
\newcommand{\scs}{\text{SCS}}
\newcommand{\mlcs}[1]{m_{\lcs(u,v),#1}}
\newcommand{\mscs}[1]{m_{\scs(u,v),#1}}
\maketitle
\begin{abstract}
A \emph{rainbow cycle} in a list-edge-colored graph is a cycle where no color appears more than once. We give a short and simple proof of a result of \cite{ABSZZ} showing that a properly edge-colored $n$-vertex graph with average degree $\Omega(\log n\log\log n)$ always contains a rainbow cycle, and generalize it to weighted graphs with edges list colored with list size at most $\ell$ and where the total weight of all edges containing a given color $i$ incident to a vertex is at most $\Delta$, showing that $\Omega_\ell(\Delta\log n\log\log n)$ average degree guarantees a rainbow cycle. 

Using the reduction of \cite{AG24} and the new generalized rainbow cycle theorem, we prove that for $r\geq 3$, every binary linear $r$-query locally correctable code $E:\F_2^k\mapsto \F_2^n$ satisfies $k\leq \tilde{O}(n^{1-2/(r-1)})$ by showing a rainbow analogue of the odd $r$-LDC lower bounds of \cite{BHKL25,JM25}, closing the gap between the best known lower bounds for odd- and even-query binary linear $r$-LCCs.
\end{abstract}
\section{Introduction}
Graph coloring is a fundamental area of study in extremal combinatorics, and one question with a large amount of recent progress is the structure of properly edge-colored graphs. A cycle in an edge-colored graph is called \emph{rainbow} if no color appears more than once, and the work of \cite{KMSV07} raised the question of the number of edges (equivalently, the average degree) needed to guarantee the existence of a rainbow cycle. Their work showed a lower bound on the average degree of $\Omega(\log n)$, and a long line of work \cite{DLS13,Janzer2020RainbowTN,Tomon2024Robust,Kim2024, JanzerSudakov} culminating in \cite{ABSZZ} proved this is nearly tight, showing that properly edge-colored graphs with average degree $\Omega(\log n\log\log n)$ must contain a rainbow cycle (closing this $\log\log n$ gap remains an open question). The rainbow cycle theorem of \cite{ABSZZ} has found many applications, notably in coding theory. 

A binary code $E:\F_2^k\rightarrow\F_2^n$ is called a \emph{$(r,\delta)$-locally correctable code}, abbreviated as $(r,\delta)$-LCC, if any codeword $E(x)$ has at most a $\delta$ fraction of coordinates flipped to produce some $y\in \F_2^n$, for each coordinate $i\in [n]$, there exists an algorithm which queries at most $r$ coordinates of $y$ and outputs the $i$th bit $E(x)_i$ of the original codeword with nontrivial probability. LCCs were introduced for program checking \cite{BK95}, and played a central role in the proof of the PCP theorem \cite{AS98,ALMSS98}, among numerous other applications. A central question in the study of LCCs is determining the optimal \emph{rate} (the ratio of the message length $k$ to the blocklength $n$ of a codeword) of a $r$-query LCC. 

Until recently, lower bounds on LCCs were simply those of locally decodable codes, or LDCs, where instead of the code admitting $r$-query decoding algorithms which can output the $i$th bit of the codeword $E(x)$, we only need to have $r$-query decoding algorithms for each bit $x_i$, $i\in [k]$ of the message. $1$-LDCs do not exist \cite{KT00} and for $2$-LDCs, we have an exponential lower bound of $k\leq \tilde O(\log n)$ on the blocklength, which is tight \cite{KdW04,GKST06}. These results carry over to LCCs. For even $r\geq 4$-LDCs, the best lower bounds are $k\leq\tilde  O(n^{1-2/r})$ \cite{Woo07,Woo08,Woo12,BCG20}, and for odd $r$-LDCs, the best known bounds were simply those of even $(r+1)$-LDCs (since any $r$-LDC is also a $(r+1)$-LDC) until a breakthrough work of \cite{AlrabiahGKM23} proved $k\leq\tilde  O(n^{1-2/r})$ for $r=3$ using the \emph{Kikuchi Method}, at the same time also giving a simple proof of the even $r$ case. The follow-up works of \cite{BHKL25,JM25} generalized the Kikuchi method to obtain the same bound for all odd $r\geq 5$. For linear $3$-LCCs, \cite{KothariM23} obtained an exponential $k\leq O(\log^8 n)$ lower bound, which was later improved to $k\leq O(\log^4 n)$ by \cite{Yankovitz24} (see also \cite{Manohar24}), which, prior to the rainbow cycle theorem, were the only results showing a stronger lower bound for linear $r$-LCCs over a constant-sized $\F_q$ compared for those known for $r$-LDCs.

The work of \cite{AG24} showed a surprising connection between LCCs and rainbow cycles, proving a reduction from binary linear $r$-LCC lower bounds to a stronger notion of "rainbow" $(r-1)$-LDC lower bounds. By using  the rainbow cycle theorem of \cite{ABSZZ}, they prove rainbow lower bounds for even $r$-LDCs, implying a near-optimal bound of $k\leq O(\log^2n\log\log n)$ for binary linear $3$-LCCs, nearly matching the $k=O(\log^2 n)$ given by Reed-Muller codes. They also obtain an improved $k\leq\tilde  O(n^{1-2/(r-1)})$ lower bound for odd $r$-LCCs. However, this does not work for even $r$-LCCs, as we cannot reduce odd LDC lower bounds to analyzing edge-colored graphs where each edge is assigned exactly one color; instead, the graphs here have edges colored by a pair of colors. This in turn motivates the question of proving a generalization of the rainbow cycle to list-colored graphs.
\subsection{Our results}
Our main result is the following theorem, which generalizes the result of \cite{ABSZZ} to weighted, locally bounded list-colored graphs. 
\begin{restatable}{theorem}{fullrainbow}
    \label{thm:fullrainbow}
Let $\cC$ be a set of colors and $G$ be an edge-weighted $N$-vertex graph, where each edge $e\in E(G)$ is assigned a nonempty list $L(e) \subseteq \cC$ where $|L(e)| \leq \ell$, and a weight $w(e)\in \R_{\geq 0}$. Suppose that, for each color $c \in \cC$ and
any vertex $u \in V(G)$, the total weight $\deg_c(u)=\sum_{e\ni S, L(e)\ni c}w(e)$ of all edges $e$ incident to $u$ such that $L(e)\ni c$ is at most $\Delta$. For sufficiently large $N$, if
$d(V(G)) \geq 13\Delta\ell 2^{\ell}\log(N)\log\log(N)$, then $G$ contains a rainbow cycle.
\end{restatable}
Using this theorem, we close the gap between even and odd $r$-LCCs, proving the following lower bound for all $r$.
\begin{restatable}{theorem}{lccbound}
    \label{thm:lccbound}
    Let $E:\F_2^k \rightarrow\F_2^n$ be a linear $(r,\delta)$-LCC for any $r\geq 3$. Then we have $k\leq O(\delta^{-2}n^{1-2/(r-1)}\log^3 n)$.
\end{restatable}
\subsection{Technical Overview}
In this section, we give an overview of the proofs of \Cref{thm:fullrainbow,thm:lccbound}. We first summarize the Kikuchi method of \cite{AlrabiahGKM23,HKMMS24,AG24,BHKL25,JM25} for LDC/LCC lower bounds and the approach of \cite{ABSZZ} for the rainbow cycle theorem.

\parhead{Kikuchi graphs and LDCs/LCCs} Given a linear $r$-LDC, standard reductions allow us to prove lower bounds for \emph{normal linear LDCs}, which have the property that for each bit $i\in [k]$ of the message, there exists a $r$-uniform size-$\Omega(n)$ hypergraph matching $\cH_i$ on $[n]$ such that to decode the $i$th bit $x_i$ of a message $x$ given the encoding $E(x)$, we simply pick any $C\in \cH_i$, and take the XOR of the coordinates in $C$ of $E(x)$. We note that if we have an even cover (a set of hyperedges $C_1, \cdots, C_t$ such that $\oplus_{u\in [t]}C_u=\emptyset$, i.e. it cover every vertex an even (possibly $0$) number of times), such that there is some $\cH_i$ such that we use an odd number of hyperedges from $\cH_i$, then taking the XOR of all $E(x)_i$, $i\in C_u$, will return $0$. However, letting $C_u\in \cH_{i_u}$ for all $u$, the normal LDC decoding property implies that this is equal to the XOR of the $i_u$th bit off $x$, and since there is some $x_i$ that appears an odd number of times, this is not equal to $0$ if $x_i=1$. Therefore, proving a bound of $k<k(r,n)$ for all $r$-LDCs is equivalent to showing that when $k\geq k(r,n)$, any set of $\Omega(n)$, $r$-uniform hypergraph matchings contains an even cover using some $\cH_i$ an odd number of times.

For even $r$, \cite{AlrabiahGKM23} constructs a level$-\ell$ Kikuchi graph, which has vertex set consisting of all size-$\ell$ subsets of $[n]$, and edges between any two $S,T\in {[n]\choose \ell}$ such that the symmetric difference $S\oplus T$ is a hyperedge. Any cycle in the Kikuchi graph corresponds to an even cover in the set of hypergraph matchings, and by coloring an edge $(S,T)$ with color $i$ if $S\oplus T\in \cH_i$, a sufficient condition for the desired even cover property is to find a cycle which uses some color exactly once, which we call an \emph{weak rainbow cycle}. By deleting a $\frac12$ fraction of edges, we can make the Kikuchi graph properly edge colored, after which classic results \cite{KMSV07} show the existence of a weak rainbow cycle. 

For odd $r$-LDCs, the aforementioned definition of Kikuchi graphs does not work: the symmetric difference of any two size-$\ell$ sets is always even and thus can never be a hyperedge in one of the odd $r$-uniform hypergraph matchings we have. Instead, \cite{AlrabiahGKM23,HKMMS24,BHKL25,JM25} consider pairs of hyperedges $C\in \cH_i,C'\in \cH_j$ from separate matchings with nonempty intersection, instead connecting two Kikuchi vertices $S,T$ if $S\oplus T=C\oplus C'$ and coloring the edge with a pair of colors $\{i,j\}$\footnote{For some, but not all, collections of hypergraph matchings \cite{JM25} may avoid this by constructing a bipartite Kikuchi graph whose edges are colored with a single color.}. As before, finding a weak rainbow cycle is sufficient for proving the desired even cover property, and this can be done using results for list-colored graphs from \cite{KMSV07}. Thus, for LDCs, we have matching lower bounds for odd $r$-LDCs.

\cite{AG24} show that if for all $k\geq k'(r-1,n)$, any set of size-$\Omega(n)$, $(r-1)$-uniform hypergraph matchings $\cH_1\cup \cdots\cup \cH_k$ contains a \emph{rainbow even cover}, meaning an even cover using at most one hyperedge from each $\cH_i$, then for every binary linear $r$-LCC, we have $k\lesssim k'(r-1,n)\log n$. When $r-1$ is even, they apply the result of \cite{ABSZZ} to the Kikuchi graph for even LDCs to show the existence of a rainbow cycle and therefore a rainbow even cover. However, the rainbow cycle theorem does not apply to list-colored graphs, so they do not obtain the same result when $r-1$ is odd.
\newcommand{\davg}{d_{\text{avg}}}

\parhead{Finding rainbow cycles: current approaches and barriers to generalization} \cite{ABSZZ} makes use of \emph{robust sublinear expanders}, which are graphs in which any subset $S$ of $n^{1-\eps}$ vertices has a size $\Omega(\eps|S|)=\Omega(|S|(1-\frac{\log |S|}{\log n} ))$ neighborhood even after $\Omega(\eps d_{\text{avg}}|S|)$ edges incident to $S$ are adversarially deleted. Given a properly edge-colored graph $G$ with color set $\cC$, they find a robust sublinear expander subgraph $H$ and obtain a rainbow cycle $H$ as follows: first, pick an arbitrary vertex $x\in V(H)$, and partition $\cC=\cC'\sqcup \cC''$ uniformly at random. Then if there is a vertex $y\in V(H)$ reachable from $x$ by both a rainbow path using only colors in $\cC'$, and a rainbow path using only colors in $\cC''$, then the union of the two paths is a rainbow cycle in $H$. By showing that the set $U(\cC')$ of all vertices reachable from $x$ via a rainbow path in $\cC'$ has size $>\frac{n+1}{2}$ with probability $>\frac12$, one can show that $U(\cC')$ and $U(\cC'')$ (defined analogously, and which has the same distribution as $U(\cC')$) intersect with positive probability and therefore one can find a rainbow cycle.

To do so, after $\cC'\subseteq \cC$ is chosen uniformly at random, they construct a sequence of sets $\cC'=\cA_0\supseteq\cA_1\supseteq\cdots\cA_T$, where $T=O(\davg)$ and each element of $\cA_i$ is included in $\cA_{i+1}$ with probability $1-\frac1T$, and thus $\cA_T$ is nonempty with high probability. We fix an arbitrary $x\in U(\cA_T)$, and undo the subsampling to show that $|U(\cA_0)|>\frac{n+1}{2}$ with probability $>\frac12$. While the sets $U(\cA_{i})$ do not necessarily grow by a significant fraction each step, in expectation, it grows every two steps. Given some $\cA_i$, they show that in expectation, either $U(\cA_{i-1})$ is much larger than $U(\cA_{i})$, or $U(\cA_{i})$ is much larger than $U(\cA_{i+1})$, and thus every two steps, we grow a nontrivial amount based on the sublinear expansion of the set of reachable vertices. For any $uv$ of color $c$ from $U(\cA_{i})$ to its complement, we call it blue if $c\notin \cA_i$ and red if $c\in \cA_i$.

Conditioned on $\cA_i$, for red $uv$, we have $u\notin U(\cA_{i+1})$ if $c\notin \cA_{i+1}$, which occurs with probability $1/T$. Thus if there are sufficiently many red edges, then the robust sublinear expansion property applied to $U(\cA_{i+1})$ shows that we can find $\Omega(\eps \davg |U(\cA_{i+1})|)$ red edges crossing the boundary, such that each vertex of $U(\cA_i)$ has $O(\davg)$ of these edges. Then in expectation, $\Omega(\eps |U(\cA_{i+1})|)$ are added when we expand from $U(\cA_{i+1})$ to $U(\cA_{i})$. Otherwise, there are many blue edges, and the robust sublinear expansion property applied to $U(\cA_i)$ shows that there exist a set of $\Omega(\eps \davg |U(\cA_{i})|)$ blue edges such that each vertex of $U(\cA_i)^c$ has $O(\davg)$ incident such edges. For blue edges, we have $v\in U(\cA_{i-1})$ if $c\in \cA_{i-1}$, and this occurs with probability $\Omega(1/T)$. Thus when expanding $U(\cA_{i})$ to $U(\cA_{i-1})$, we add $\Omega(\eps |U(\cA_i)|)\geq\Omega(\eps |U(\cA_{i+1})|)$ total vertices and therefore in expectation we grow by a $1+\Omega(\eps)$ factor every two steps, so after $T=\Omega(\log n\log\log n)$ steps, we obtain a size $>\frac{n+1}{2}$ set. \cite{ABSZZ} then uses an intricate probabilistic argument to show $\Pr[|U(\cA_0)|>\frac{n+1}{2}]>\frac12$.

However, the fact that each edge is colored with exactly one color is crucial, and this breaks down when each edge is colored with two colors. In this case, we call an edge $uv$ blue if $c,c'\notin \cA_i$ and red otherwise. Here, we get $v\in U(\cA_{i-1})$ through $uv$ with probability $O(1/T^2)$, and thus we do not get $\Omega(\eps|U(\cA_i)|)$ expansion. Since there are graphs where with high probability, most edges are blue for a constant of $\cA_i$s (the hypercube, for example, where each edge is given colors $\{\alpha_i,\beta_i\}$ if the $i$th bit flips, where these colors are distinct for different $i$), this is unavoidable, and thus we cannot iteratively show expansion in this way.

\parhead{Our approach} Rather than showing that the set of reachable vertices expands as we pairs of colors, we instead pass to a sublinear expander subgraph with sufficiently large average degree and show that, after again fixing an arbitrary vertex $x$ in a sublinear expander subgraph, the expected number of vertices reachable from $x$ via a random subset $\cC_{q}$ of colors, where each color is included independently with probability $q$, is larger than $\frac{n+1}{2}$ when $q=\frac12$. We obtain an upper bound on the expected isoperimetric ratio $\frac{|\partial U(\cC_q)|}{|U(\cC_q)|}$ by bounding the size of the boundary of $U(\cC_q)$ by the sum over all colors $c$ of the change in $|U(\cC_q)|$ when a single color is added or deleted. The number of red edges is bounded by the change when some $c\in U(\cC_q)$ is deleted, while the number of blue edges is bounded by the change when some $c\notin U(\cC_q)$ is added. The Russo-Margulis Formula then upper bounds this as a function of $\frac{d}{dq}\E \log |U(\cC_q)|$. Since we have a lower bound on $\frac{|\partial U(\cC_q)|}{|U(\cC_q)|}$ from sublinear expansion, a large average degree implies that $\E|U(\cC_q)|$ grows quickly as a function of $q$ and thus $\E|U(\cC_{1/2})|$ is large. 

This approach generalizes to the list-coloring setting. With list size $\ell$, red edges can be bounded by the change in $|U(\cC_q)|$ when at most $\ell$ colors are deleted. We fix an arbitrary order of removal and consider the pivotal color, and bound by the change when that color is deleted. Similarly, for blue edges, we fix an arbitrary order of added colors and bounded by the change when the pivotal color is deleted. Taking expectations, we again can obtain an upper bound in terms of the sum of the changes when adding or deleting a single $c\in U(\cC_q)$, and obtain the same conclusion as before.

\section{A Short Proof of the Rainbow Cycle Theorem}
In this section, we give a short proof of Theorem~1.1 of \cite{ABSZZ}. Here we will present a short proof of \cite[Theorem~1.1]{ABSZZ} using our new proof approach.
\begin{theorem}\label{thm:shortabszz}
    Any properly edge-colored $N$-vertex graph $G$ with average degree $d(V(G))\geq 9\log N\log\log N$ contains a rainbow cycle.
\end{theorem}
We will make use of the following lemma to find a sublinear expander, which we prove in \Cref{sec:rainbowproofs}.
\begin{restatable}{lemma}{maxsubgraph}\label{lem:maxsubgraph}
        Let $G$ be a graph and let $H$ be a subgraph of $G$ maximizing $\frac{d(H)}{\log|V(H)|}$. Then we have $\frac{|\partial_H S|}{|S|}\geq\frac{d(H)}{2\log n}\log\frac{n}{|S|}$ for all $\emptyset \neq S\subseteq V(H)$.

\end{restatable}
\begin{proof}[Proof of \Cref{thm:shortabszz}]
    Let $\cC$ be the set of all colors and let $H$ be a subgraph of $G$ maximizing $\frac{d(H)}{\log|V(H)|}$ and let $n\coloneqq |V(H)|$. Fix an arbitrary vertex $x\in V(H)$ and for any $\cA\subseteq \cC$, let $U(\cA)=U_x(\cA)$ denote the set of vertices in $H$ reachable from $x$ by a rainbow walk using colors from $\cA$. 
    
    Let $0\leq q\leq \frac12$ and $\cC_q$ be a random subset of $\cC$ such that every $c\in \cC$ is included in $\cC_q$ independently with probability $q$. Define $f(q)=\E\log\frac{n}{|U(\cC_q)|}=\log n-\E \log|U(\cC_q)|$.

    Call an edge $uv$ of color $c$ from $U(\cC_q)$ to its complement \emph{red} if $u\notin U(\cC_q\setminus c)$, i.e. every path from $x$ to $u$ uses color $c$, and \emph{blue} otherwise, i.e. one can reach $u$ without $c$. For every red $uv$ of color $c$, we have $u\in U(\cC_q)\setminus U(\cC_q\setminus c)$, and since $H$ is properly edge colored, the total number of red $uv$ is at most $\sum_{c\in \cC_q}|U(\cC_q)|-|U(\cC_q)\setminus c|$. On the other hand, for every blue $uv$ of color $c$, we have $u\in U(\cC_q)$ and $v\in U(\cC_q\cup c)\setminus U(\cC_q)$, so the number blue edges of color $c$ is at most $\min(|U(\cC_q)|,|U(\cC_q\cup c)|-|U(\cC_q)|)$ and thus there are at most $\sum_{c\notin \cC_q}\min(|U(\cC_q)|,|U(\cC_q\cup c)|-|U(\cC_q)|)$ blue edges in total. Therefore, we have
    \begin{align*}
        \E \frac{|\partial U(\cC_q)|}{|U(\cC_q)|}
        &\leq \E \sum_{c\in \cC_q}\frac{|U(\cC_q)|-|U(\cC_q)\setminus c|}{|U(\cC_q)|}+\E\sum_{c\notin \cC_q}\frac{\min(|U(\cC_q)|,|U(\cC_q\cup c)|-|U(\cC_q)|)}{|U(\cC_q)|}
        \\&=\E \sum_{c\in \cC_q}\left(1-\frac{|U(\cC_q)\setminus c|}{|U(\cC_q)|}\right)+\E\sum_{c\notin \cC_q}\min\left(1, \frac{|U(\cC_q\cup c)|}{|U(\cC_q)|}-1\right)
        \\&\leq \E \sum_{c\in \cC_q}\log\frac{|U(\cC_q)|}{|U(\cC_q\setminus c)|} + 2\E\sum_{c\notin \cC_q}\log \frac{|U(\cC_q\cup c)|}{|U(\cC_q)|}
        \\&\leq \E\sum_{c\in \cC_q}\log |U(\cC_q)|-\log|U(\cC_q\setminus c)| + 2\E\sum_{c\notin \cC_q}\log |U(\cC_q\cup c)|-\log|U(\cC_q)|
        \\&\leq q\E \sum_{c\in \cC}\log |U(\cC_q)\cup c|-\log |U(\cC_q\setminus c)| + 2(1-q)\E \sum_{c\in \cC}\log |U(\cC_q)\cup c|-\log |U(\cC_q\setminus c)|
        \\&\leq 2\frac{d}{dq}\E \log |U(\cC_q)|=-2\frac{d}{dq}f(q)\mcom
    \end{align*}
    where we use $\min(1,z-1)\leq 2\log z$ in the third line and the last inequality follows from the Russo-Margulis Formula (\Cref{fact:russomargulis}).
    
    By \Cref{lem:maxsubgraph}, we have $\E\frac{|\partial U(\cC_q)|}{|U(\cC_q)|}\geq \E \frac{d(H)}{2\log n}(\log n-\log |U(\cC_q)|)= \frac{d(H)}{2\log n}f(q)$. Therefore
    \begin{align*}
         f'(q)=-\frac{d}{dq}\E\log |U(\cC_q)|
         \leq -\frac12\E\frac{|\partial U(\cC_q)|}{|U(\cC_q)|}
         \leq -\frac{d(H)}{ 4\log n}f(q)\mcom
    \end{align*}
    which implies that $f(\frac12)\leq f(0)\exp\left(\frac{-d(H)}{8\log n}\right)\leq \exp\left(-\frac{ d(H)}{8\log n}\right)\log n$.

    Suppose $G$ does not contain a rainbow cycle, in which case $U(\cC_q)\cap U(\cC_q^c)=\{x\}$, so $|U(\cC_q)|+|U(\cC_q^c)|\leq n+1$ and $\log|U(\cC_q)|+\log|U(\cC_q^c)|= \log(|U(\cC_q)||U(\cC_q^c)|)\leq  2\log((n+1)/2)$. Since $\cC_{1/2}\overset{d}{=}\cC_{1/2}^c$, this implies that $\E\log U(\cC_{1/2})\leq \log(\frac{n+1}{2})$, so $f(\frac12)\geq \log\frac{2n}{n+1}\geq \log(4/3)=\Omega(1)$. Combined with above, we get
    \begin{align*}
        d(H)
        \leq 8\log n\log\log \left(\frac{n}{f(1/2)}\right) &\leq(8+o(1))\cdot\log n\log\log n
        < 9\log n\log\log n\mcom
    \end{align*}
    and therefore $d(G)\leq \frac{\log N}{\log n}d(H) <9\log N\log\log n\leq 9\log N\log\log N.$
\end{proof}

\section{Preliminaries}
\newcommand{\dec}{\mathsf{Dec}}
\newcommand{\krain}[1]{k_{\text{rainbow}}^{(#1)}}
We first formally define a $r$-query locally correctable code.
\begin{definition}[Binary Linear LCC] A linear code $E : \F_2^k\rightarrow \F_2^n$ is called a $(r,\delta)$-locally correctable code (LCC) if there exists a randomized decoding algorithm $\dec^y(\cdot)$ with oracle access to $y\in \F_2^n$, takes in an input $i\in [n]$, and satisfies the following properties:
\begin{itemize}
    \item $\dec^y(i)$ makes at most $r$ queries to string $y$.
    \item For all $b\in \F_2^k$ and for all $y\in \{\pm 1\}^n$ such that $\operatorname{dist}(y, E(b))\leq \delta n$, $\Pr[\dec(i)=E(b)_i]\geq \frac23$
\end{itemize}
\end{definition}
\begin{definition}[Normal Linear LCC]
    Let $E : \{\pm 1\}^k\rightarrow \{\pm 1\}^n$ be a linear code with generator matrix $M\in \F_2^{n\times k}$ and let $v_i$ denote the $i$th row of $M$ for all $i\in [n]$. We call $E$ a normal $(r,\delta)$-locally correctable code (LCC) if there exists $r$-uniform hypergraph matchings $\cH_1, \cdots, \cH_n$ over vertex set $[n]$ such that $|\cH_i|\geq \delta n$ for all $i\in [n]$ and for all $i\in [n]$ $C=\{a_1, \cdots, a_r\}\in \cH_i$, we have $v_i=\sum_{j=1}^r v_{a_j}$.
\end{definition}
We have the following reduction from linear LCCs to normal LCCs.
\begin{fact}
\label{fact:normallccreduction}
    Let $E:\F_2^k\rightarrow \F_2^n$ be a linear $(r,\delta, \eps)$-locally correctable code. Then there exists $E':\F_2^k\rightarrow \F_2^{2n}$ such that $E'$ is a linear normal $(r, \frac{\delta}{2r})$-locally correctable code.
\end{fact}

\begin{definition}
    Let $\delta > 0$ and $q,n\in \N$ such that $q\geq 2$. We define $\krain{r}(\delta, n)$ to be the smallest $k$ such that if for any collection $\cH_1, \cdots, \cH_k$ of $r$-matchings over vertex set $[n]$ such that $|\cH_i|\geq \delta n$ for all $i\in [n]$, there exists some $\emptyset\neq \mathcal{E}\subseteq \bigcup_{i=1}^k\cH_i$ such that $\bigoplus_{C\in \cE} C=\emptyset$ and $|\cE\cap \cH_i|\leq 1$ for all $i\in [k]$.
\end{definition}
\cite{AG24} proves the following upper bound on $\krain{q}$ for even $q$.
\begin{fact}[Proposition~4.3 of \cite{AG24}]\label{fact:evenrainbowbound}
    For any even $q\geq 4$ and $\delta\in (0,1)$, we have $\krain{q}(\delta,n)\leq O(\delta^{-1}n^{1-2/q}\log^2 n)$
\end{fact}

We will use the following result of \cite{AG24} to prove lower bounds on even-query locally correctable codes.
\begin{fact}[Theorem~4.2 of \cite{AG24}]\label{fact:rainbowtolcc} Let $E:\{\pm 1\}^k \rightarrow\{\pm 1\}^n$ be a normal linear $(r,\delta)$-LCC for any $r\geq 3$. Then 
$k\leq O(\delta^{-1}\log n \cdot k_{\text{rainbow}}^{r-1}(\delta/2,n)).$    
\end{fact}
\begin{fact}[Lemma~2.7 of \cite{HKMMS24}] 
\label{fact:even-cover}
    Suppose we have a normal $q$-LDC with $q$-uniform matchings $\cH_1, \cH_2, \dots, \cH_k$ over $[n]$ and we color each hyperedge in $\cH_i$ with color $i$. 
    Then, any even cover in $\cH = \bigcup_{i\in [k]} \cH_i$ must use color $i$ an even number of times for each $i\in [k]$.
\end{fact}
In addition, we use the following result.
\begin{fact}[Russo-Margulis Formula, see Theorem~2.32 of \cite{Grimmett}]\label{fact:russomargulis}
    Let $\cS$ be a finite set and $\cS_p$ be the random subset of $\cS$ obtained by including every element independently with probability $p$. For every function $f:2^{\cS}\rightarrow \R$, we have
    \begin{align*}
        \frac{d}{dq} \E f(\cS_q)= \E\sum_{s\in \cS}f(\cS_q\cup c)-f(\cS_q\setminus c)\mper
    \end{align*}
\end{fact}
\subsection{Kikuchi Graphs for odd $r$-LDCs}
Here we define the Kikuchi graphs used in \cite{BHKL25}.
\begin{definition}
    For any hypergraph $\cH$ with vertex set $[n]$, the co-degree of a set $Q\subseteq [n]$, which we denote as $d_{\cH,Q}$, is the number of hyperedges in $\cH$ containing $Q$. Moreover, we denote $d_{\cH,t} \coloneqq \max_{|Q| =t} d_{\cH,Q}$ for all $t\in [q]$. When the context is clear, we omit the dependency on $\cH$ and write $d_t$ for simplicity.
\end{definition}

\begin{definition}[Good index $t$, Definition 2.6 of \cite{BHKL25}]
\torestate{
\label{def:goodindex}
Let $d_1\geq d_2 \geq \cdots \geq d_q$. We say that an index $t$ is \emph{good} with respect to the tuple $(d_1, d_2, \ldots, d_q)$ if the following conditions hold:
\begin{enumerate}[(1)]
\item $d_r/d_{t}\leq n^{1-\frac{2r}{q}}$ for every $1 \leq r \leq \lceil \frac{q-t}{2} \rceil$, \label{item:good1}
\item $d_r/d_{t} \leq n^{-\frac{2}{q}(r-t) + \frac{1}{q}(t-\1(t\text{ even}))}$ for every $t \leq r \leq \lfloor \frac{q+t}{2} \rfloor$, \label{item:good2}
\item If $t <q/2$, then $d_{t} \geq d_1 n^{-\frac{2}{q}(t-1)}$, and if $t>q/2$, then $d_{t} \geq d_1 n^{-1+2/q}$. \label{item:d-t-lower-bound} \label{item:good3}
\end{enumerate}
}
\end{definition}
\begin{definition}[Approximate Strong Regularity, Definition 2.7 of \cite{BHKL25}]
\label{def:approxstrongregdef}
For a given $t\in[q]$, a $q$-uniform hypergraph $\cH$ on $[n]$ is called $t$-approximately strongly regular if there exists a partition $\cH = \cH^{(t)}_1 \sqcup \cdots \sqcup \cH^{(t)}_{p_t}$ such that:
    \begin{enumerate}[(1)]
        \item For every $\theta\in[p_t]$, there is some $Q_\theta \in \binom{[n]}{t}$ such that all $C \in \cH^{(t)}_\theta$ contain $Q_\theta$.
        \item $d_{t}/2\leq \big|\cH^{(t)}_\theta\big|\leq d_{t}$ for all $\theta\in[p_t]$.
        \item $t$ is a good index 
        with respect to the sequence $d_1 \geq \cdots \geq d_q$ (recall that $d_r:= \max_{|Q| = r} d_{\cH,Q}$). \label{item:good-index}
    \end{enumerate}

\end{definition}
\newcommand{\chtt}{\cH_{\theta}^{(t)}}
\newcommand{\chtp}{\cH_{\theta}^{\pi}}
\begin{definition}[Definition 2.11 of \cite{BHKL25}] \label{def:kikuchi}
Given a $t$-approximately strongly regular hypergraph $\cH$ on vertex set $[n]$, hypergraph matchings $\cH = (\cH_1,\dots,\cH_k)$, and a partition $\{\cH_{\theta}^{(t)}\}_{\theta\in [p_t]}$ from approximate strong regularity, we define the level-$\ell$ Kikuchi graph $K_\ell(\cH)$ as follows. 

We create two copies $(u,1)$ and $(u,2)$ of the original $n$ variables $1 \leq u \leq n$. The Kikuchi graph has vertices indexed by subsets $S \in {[n] \times [2] \choose \ell}$ for a parameter $\ell \in \N$.
For any $C \neq C' \in \cH^{(t)}_{\theta}$, let  $\widetilde{C} \coloneqq C\setminus Q_{\theta}$ and $\widetilde{C}' \coloneqq C'\setminus Q_{\theta}$.
For any $S,T \in {[n] \times [2]\choose \ell}$, we draw an edge $S\xleftrightarrow{C, C'} T$ with colors $\gamma(S,T)=(i,j)$ if:
    \begin{enumerate}
        \item $S\oplus T = X\oplus Y$, for $X \coloneqq \widetilde{C}\times\{1\}, Y \coloneqq \widetilde{C}'\times\{2\}$. Clearly, $X\oplus Y = X\cup Y$ as $X, Y$ are disjoint.
        \item When $q - t$ is even, $|S\cap X| = |S\cap Y| = |T\cap X| = |T\cap Y| = (q - t)/2$.
    \item When $q - t$ is odd, either
    \begin{enumerate}
        \item $|S\cap X| = |T\cap Y| = \left\lceil\frac{q - t}{2}\right\rceil, |S\cap Y| = |T\cap X| = \left\lfloor\frac{q - t}{2}\right\rfloor$ or
        \item $|S\cap X| = |T\cap Y| = \left\lfloor\frac{q - t}{2}\right\rfloor, |S\cap Y| = |T\cap X| = \left\lceil\frac{q - t}{2}\right\rceil$.
    \end{enumerate}
    \end{enumerate}
    Let $\deg^{\cH}_i(S)=\#\{(S,T)\in E(K_\ell(\cH)): \gamma(S,T)\ni i\}$ and $\deg^\cH(S)=\sum_{i\in [k]}\deg^{\cH}_i(S)$.
\end{definition}
\section{Proofs of \Cref{thm:fullrainbow} and \Cref{thm:lccbound}}
In this section we prove \Cref{thm:fullrainbow} and \Cref{thm:lccbound}. We begin with \Cref{thm:fullrainbow}, which we restate below.
\fullrainbow*
We now prove the theorem using \Cref{lem:maxsubgraph}, which we restate again, and our key technical lemma using the Russo-Margulis Formula, both of which we prove in \Cref{sec:rainbowproofs}.
\begin{lemma}\label{lem:maxsubgraph}
    Let $G$ be a graph and let $H$ be a subgraph of $G$ maximizing $\frac{d(H)}{\log|V(H)|}$. Then we have $\frac{|\partial_H S|}{|S|}\geq\frac{d(H)}{2\log n}\log\frac{n}{|S|}$ for all $\emptyset \neq S\subseteq V(H)$.
\end{lemma}
\begin{lemma}\label{lem:expansion}
    Consider the setting of \Cref{thm:fullrainbow} and let $\cC_q\subseteq \cC$ be a random subset where each $c\in \cC$ is included independently with probability $q$. For an arbitrary vertex $x$, let $U(\cC_q)$ denote the set of all vertices reachable by a rainbow walk from $x$. Then we have 
\begin{align*}
    \E \frac{|\partial U(\cC_q)|}{|U(\cC_q)|}\leq \Delta\left(\frac{q}{(1-q)^{\ell-1}}+\frac{2(1-q)}{q^{\ell-1}}\right)\frac{d}{dq}\E \log|U(\cC_q)|\mper
\end{align*}
\end{lemma}
\begin{proof}[Proof of \Cref{thm:fullrainbow}]
        \newcommand{\nbc}{\nu_{\cB,c}}
    Let $H$ be a subgraph of $G$ maximizing $\frac{d(H)}{\log|V(H)|}$ and let $n\coloneqq |V(H)|$. Fix an arbitrary vertex $x\in V(H)$ and for any $\cA\subseteq \cC$, let $U(\cA)=U_x(\cA)$ denote the set of vertices in $H$ reachable from $x$ by a rainbow walk using colors from $\cA$. Define $h(U(\cA))\coloneqq \frac{|\partial U(\cA)|}{|U(\cA)|}$. 
    
    Let $\frac12-\frac{1}{6\ell}\leq q\leq \frac12$ and $\cC_q$ be a random subset of $\cC$ such that every $c\in \cC$ is included in $\cC^q$ independently with probability $q$. Define $f(q)=\log n-\E \log|U(\cC_q)|$. By \Cref{lem:maxsubgraph}, we have $\E\frac{|\partial U(\cC_q)|}{|U(\cC_q)|}\geq \E \frac{d(H)}{2\log n}(\log n-\log |U(\cC_q)|)= \frac{d(H)}{2\log n}f(q)$, while by \Cref{lem:expansion}, we have $\E \frac{|\partial U(\cC_q)|}{|U(\cC_q)|}\leq \Delta\left(\frac{q}{(1-q)^{\ell-1}}+\frac{2(1-q)}{q^{\ell-1}}\right)\frac{d}{dq}\E \log|U(\cC_q)|\leq \Delta 2^{\ell}\frac{d}{dq}\E\log|U(\cC_q)|$. Therefore
    \begin{align*}
         f'(q)=-\frac{d}{dq}\E\log |U(\cC_q)|
         \leq -\frac{1}{\Delta 2^\ell}\E\frac{|\partial U(\cC_q)|}{|U(\cC_q)|}
         \leq -\frac{d(H)}{\Delta 2^{\ell+1}\log n}f(q)\mcom
    \end{align*}
    which implies that $f(\frac12)\leq f(\frac12-\frac{1}{6\ell})\exp\left(\frac{-d(H)}{6\ell\Delta2^{\ell+1}\log n}\right)\leq \exp\left(-\frac{ d(H)}{6\ell\Delta 2^{\ell+1}\log n}\right)\log n$.

    Suppose $G$ does not contain a rainbow cycle, in which case $U(\cC_q)\cap U(\cC_q^c)=\{x\}$, so $|U(\cC_q)|+|U(\cC_q^c)|\leq n+1$ and $\log|U(\cC_q)|+\log|U(\cC_q^c)|= \log(|U(\cC_q)||U(\cC_q^c)|)\leq  2\log((n+1)/2)$. Since $\cC_{1/2}\overset{d}{=}\cC_{1/2}^c$, this implies that $\E\log U(\cC_{1/2})\leq \log(\frac{n+1}{2})$, so $f(\frac12)\geq \log\frac{2n}{n+1}\geq \log(4/3)=\Omega(1)$. Combined with above, we get
    \begin{align*}
        d(H)
        \leq 6\ell\Delta 2^{\ell+1}\log n\log\log \left(\frac{n}{f(1/2)}\right) &\leq(1+o(1))6\ell\Delta 2^{\ell+1}\cdot\log n\log\log n
        \\&\leq 13\ell\Delta 2^\ell\log n\log\log n\mper
    \end{align*}
    Therefore, $d(G)\leq \frac{\log N}{\log n}d(H)\leq 13\ell\Delta 2^{\ell}\log N\log\log N$.
\end{proof}
We then show \Cref{thm:lccbound} given the following lemma whose proof is deferred to \Cref{sec:lemproof}.
\begin{lemma}\label{lem:rainbowbound}
    For any odd $q\geq 3$ and $\delta\in (0,1)$, we have $k_{\text{rainbow}}^{q}(\delta,n))\leq O(\delta^{-1}n^{1-2/q}\log^2 n)$.
\end{lemma}
\lccbound*
\begin{proof}[Proof of \Cref{thm:lccbound}]
    By \Cref{fact:normallccreduction}, there exists a normal linear $(r,\frac{\delta}{2r})$-LCC $E':\F_2^k\rightarrow \F_2^{2n}$. Thus \Cref{lem:rainbowbound} (when $r$ even) or \Cref{fact:evenrainbowbound} (when $r$ odd) and \Cref{fact:rainbowtolcc} implies that \[k\leq O(\delta^{-1}\log n\cdot \krain{r-1}(\delta/4r,n))\leq O(r\delta^{-2}n^{1-2/(r-1)})\log^3 n\mper\]
\end{proof}
\subsection{Proofs of \Cref{lem:maxsubgraph,lem:expansion}}
\label{sec:rainbowproofs}
\begin{proof}[Proof of \Cref{lem:maxsubgraph}]
    For all $\emptyset\neq S\subseteq |V(H)|$, we have $d(H(S))\leq d(H)\frac{\log|S|}{\log n}$ and $d(H(\bar{S}))\leq d(H)$. Therefore,
    \begin{align*}
        |\partial S|&=|E(H)|-|E(S)|-|E(V(H)\setminus S)|
        \\&\geq \frac{nd(H)}{2}-|S|d(H)\frac{\log|S|}{2\log n}-\frac12(n-|S|)d(H)
        \\&\geq \frac{|S|d(H)}{2\log n}(\log n-\log|S|)\mper
    \end{align*}
\end{proof}
\begin{proof}[Proof of \Cref{lem:expansion}]
    For any $\cA\subseteq \cC$, call an edge $uv$ from $U(\cA)$ to its complement \emph{red} if $u\notin U(\cA\setminus L(uv))$, i.e. every path from $x$ to $u$ uses at least one element of $L(uv)$, and \emph{blue} otherwise. Define $h_r(U(\cA))$ and $h_b(U(\cA))$ to be the contributions to $h(U(\cA))$ from red and blue elements of $\partial U(\cA)$, respectively, so $h(U(\cA))=h_r(U(\cA))+h_b(U(\cA))$. We will upper bound the quantity $\E h(U(\cC^q))$ by considering red and blue edges separately. 

    Let $uv$ be a red edge, so $u\notin U(\cC^q\setminus L(uv))$. Then there exists some $\cC^q\setminus L(uv)\subset \cA\subseteq \cC^q$ and $c\in \cA\cap L(uv)$ such that $u\in U(\cA)\setminus U(\cA\setminus c)$. Letting $\cB=\cC^q\setminus \cA$, this is equivalent to having some $c\in L(uv)$, $\cB\subset L(uv)\setminus c$, such that $u\in U(\cC^q\setminus \cB)\setminus U(\cC^q\setminus(\cB\cup c))$. Assign each red edge to one such pair $\cA,c$. Note that for each $v\in U(\cA)\setminus U(\cA\setminus c)$, the total weight of all such edges $e$ incident to $v$ such that $c\in L(e)$ is at most $\Delta$, so the total weight of red edges mapped to $\cA,c$ is at most $\Delta(|U(\cA)|-|U(\cA\setminus c)|)$, so the total weight of all edges assigned to $\cA$ is at most $\Delta\sum_{c\in \cA}|U(\cA)|-|U(\cA\setminus c)|$.

    Clearly, any such $\cA$ must satisfy $|\cA|\geq |\cC_q|-\ell+1$. To bound $\E h_r(\cC_q)$, first condition on $|\cC_q|=T$, in which case picking $\cC_q$ followed by $\cA$
    is equivalent to first picking $0\leq a\leq \ell-1$ and $\cA$ such that $|\cA|=T-a$, then choosing the $a$ elements of $\cC_q\setminus \cA\coloneqq \cB$. Letting $\eta_{\cB}(\cA)$ denote the total weight of all edges $uv\in \partial U(\cA\cup \cB)$ such that there exists some $c\in \cA\cap L(uv)$ with $\cB\subset L(uv)\setminus c$ and $u\in U(\cA)\setminus U(\cA\setminus c)$, note that for each $uv$ there are at most ${|L(uv)\setminus c|\choose a}\leq {\ell-1\choose a}$ total $\cB\in \cC\setminus \cA$ such that $uv$ contributes to $\eta_\cB(\cA)$. Therefore we have for all $a,\cA$,
    $$\sum_{\substack{\cB\in \cC\setminus\cA\\|\cB|=a}}\eta_{\cB}(\cA)\leq \Delta{\ell - 1\choose a}\sum_{c\in \cA}|U(\cA)|-|U(\cA\setminus c)|\mper$$

    Finally, note that $\Pr(\cC_q=\cA\cup \cB)=\left(\frac{q}{1-q}\right)^a\Pr(\cC_q=\cA)$ for all $|\cB|=a$ independent of $\cA$. Therefore,

    \begin{align*}
        \E h_r(U(\cC_q))
        &\leq
        \sum_{\cA}\sum_{a=0}^{\ell-1}\sum_{\substack{\cB\in \cC\setminus\cA\\|\cB|=a}}\Pr(\cC_q=\cA\cup \cB)\cdot \frac{\eta_\cB(\cA)}{|U(\cA\cup \cB)|}
        \\&\leq \Delta\sum_{\cA}\sum_{a=0}^{\ell-1}\left(\frac{q}{1-q}\right)^a{\ell-1\choose a}\Pr(\cC_q=\cA)\sum_{c\in \cA}\frac{|U(\cA)|-|U(\cA\setminus c)|}{|U(\cA)|}
        \\&\leq\frac{\Delta}{(1-q)^{\ell - 1}}\sum_{\cA}\Pr(\cC_q=\cA)\sum_{c\in \cA}\frac{|U(\cA)|-|U(\cA\setminus c)|}{|U(\cA)|}
        \\&\leq \frac{\Delta}{(1-q)^{\ell - 1}}\E\sum_{c\in\cC_q}\frac{|U(\cC_q)|-|U(\cC_q\setminus c)|}{|U(\cC_q)|}
        \\&\leq \frac{\Delta}{(1-q)^{\ell - 1}}\E\sum_{c\in\cC_q}\log |U(\cC_q)|-\log |U(\cC_q\setminus c)|
        \\&\leq \frac{q\Delta}{(1-q)^{\ell - 1}}\frac{d}{dq}\E\log |U(\cC_q)|\mcom
    \end{align*}
    where the last line follows from \Cref{fact:russomargulis}.

    Now we consider blue edges $uv$, so $u\in U(\cC^q\setminus L(uv))$ and $v\in U(\cC^q\cup L(uv))\setminus U(\cC^q)$. Then there exists some $\cC^q\subseteq \cA\subset \cC^q\cup L(uv)$ and $c\in L(uv)\setminus \cA$ such that $v\in U(\cA\cup c)\setminus U(\cA)$. First, note that for any such $uv$,
    \begin{align*}
        \frac{w(uv)}{|U(\cC_q)|}
        \leq \sum_{\substack{\cC^q\subseteq \cA\subset \cC_q\cup L(uv)\\c\in L(uv)\setminus \cA \\ v\notin U(\cA) }}w(uv)\left(\frac{1}{|U(\cA)|}-\frac{\1\{v\notin U(\cA\cup c)\}}{|U(\cA\cup c)|}\right)\mper
    \end{align*}
    
    We first condition on $|\cC_q|=T$, in which picking $\cC_q$ followed by $\cA$ is equivalent to first picking $0\leq a\leq \ell - 1$ and $\cA$ such that $|\cA|=T+a$, then choosing the size $a$ set $\cB\coloneqq\cA\setminus \cC_q $. For all $C\subset \cC$, arbitrarily assign each blue edge in $\partial U(C)$ to one $\cA'\supseteq C$ such that there exists some $c'\in L(uv)\setminus \cA'$ such that $v\in U(\cA'\cup c')\setminus U(\cA')$. Note that for each $uv$ there are at most ${|L(uv)\setminus c|\choose a}\leq {\ell - 1\choose a}$ total $\cB\subseteq L(uv)\setminus c$ and the total weight of all $uv\in \partial U(\cA)$ such that $L(uv)\ni c$ is at most $\Delta|U(\cA)|$. Thus we have
\begin{align*}
    &\sum_{c\notin \cA}\sum_{\substack{\cB\in \cA\\|\cB|=a}}\sum_{\substack{uv\in \partial U(\cA)\\\cB\cup c\in L(uv)}}w(uv)\left(\frac{1}{|U(\cA)|}-\frac{\1\{v\notin U(\cA\cup c)\}}{|U(\cA\cup c)|}\right)
        \\&\leq {\ell - 1\choose a}\Delta 
        \sum_{c\notin \cA}\left(|U(\cA)|\left(\frac{1}{|U(\cA)|}-\frac{1}{|U(\cA)\cup c|}\right)+\frac{|U(\cA\cup c)|-|U(\cA)|}{|U(\cA\cup c)|}\right)
        \\&=2{\ell - 1\choose a}\Delta \sum_{c\notin \cA}\frac{|U(\cA\cup c)|-|U(\cA)|}{|U(\cA\cup c)|}\mper
\end{align*}
Finally, note that $\Pr(\cC_q=\cA\setminus \cB)=\left(\frac{1-q}{q}\right)^a\Pr(\cC_q=\cA)$ for all $|\cB|=a$ independent of $\cA$. Letting $\nu_{\cB}(\cA)$ denote the total weight of all edges assigned to $\cA$ for $C=\cA\setminus \cB$,
\begin{align*}
        \E h_b(U(\cC_q))
        &\leq\sum_{\cA}\sum_{a=0}^{\ell -1}\sum_{\substack{\cB\subset \cA\\|\cB|=a}}\Pr(\cC_q=\cA\setminus \cB)\frac{\nu_{\cB}(\cA)}{|U(\cA\setminus \cB)|}
        \\&\leq \sum_{\cA}\Pr(\cC_q=\cA)\sum_{a=0}^{\ell-1}\left(\frac{1-q}{q}\right)^{a}\sum_{c\notin \cA}\sum_{\substack{\cB\subset \cA\\|\cB|=a}}\sum_{\substack{uv\in \partial U(\cA)\\\cB\cup c\in L(uv)}}w(uv)\left(\frac{1}{|U(\cA)|}-\frac{\1\{v\notin U(\cA\cup c)\}}{|U(\cA\cup c)|}\right)
        \\&\leq 2\Delta\sum_{\cA}\Pr(\cC_q=\cA)\sum_{a=0}^{\ell - 1}\left(\frac{1-q}{q}\right)^{a} {\ell - 1\choose a}\sum_{c\notin \cA}\frac{|U(\cA\cup c)|-|U(\cA)|}{|U(\cA\cup c)|}
        \\&\leq \frac{2\Delta}{q^{\ell-1}} \sum_{\cA}\Pr(\cC_q=\cA)\sum_{c\notin \cA}\frac{|U(\cA\cup c)|-|U(\cA)|}{|U(\cA\cup c)|}
        \\&\leq \frac{2\Delta}{q^{\ell-1}} \E\sum_{c\notin \cC_q}\frac{|U(\cC_q\cup c)|-|U(\cC_q)|}{|U(\cC_q\cup c)|}
        \\&\leq \frac{2\Delta}{q^{\ell-1}} \E\sum_{c\notin \cC_q}\log |U(\cC_q\cup c)|-\log |U(\cC_q)|
        \\&\leq \frac{2\Delta(1-q)}{q^{\ell-1}}\frac{d}{dq}\E \log|U(\cC_q)|\mcom
    \end{align*}
    where again the last line follows from \Cref{fact:russomargulis}.
Combining the two inequalities gives
\begin{align*}
    \E h(U(\cC_q))
    =\E h_r(U(\cC_q))+\E h_b(U(\cC_q))
    &\leq \Delta\left(\frac{q}{(1-q)^{\ell-1}}+\frac{2(1-q)}{q^{\ell-1}}\right)\frac{d}{dq}\E \log|U(\cC_q)|\mper
\end{align*}
\end{proof}
\subsection{Weighted Kikuchi Graphs and Proof of \Cref{lem:rainbowbound}}
\label{sec:lemproof}

We now prove \Cref{lem:rainbowbound}. We construct a \emph{weighted Kikuchi graph} from $\cH$ such that a rainbow cycle in the Kikuchi graph corresponds to a rainbow even cover in $\cH$, then show that for an appropriate choice of parameters it contains a rainbow cycle. We begin by introducing the notion of a weak regularity subpartition, which is a collection of disjoint $t$-approximate strongly regular hypergraphs (see \Cref{def:approxstrongregdef}) whose $t$-codegree decays geometrically and whose union contains a constant fraction of $|\cH|$.
\begin{definition}[Weak Regularity]\label{def:weakregdef}
    Given $t\in [q]$, a $q$-uniform hypergraph $\cH$ on $[n]$, we call a collection $\cP=\{\cH^\pi :\pi\in P\}$ of disjoint subhypergraphs of $\cH$ a \emph{$t$-weakly regular subpartition of $\cH$} if:
    \begin{enumerate}
        \item $\cH^\pi$ is $t$-approximately strongly regular for all $\pi\in P$,
        \item $d_t^\pi\geq 3$ for all $\pi\in P$,
        \item There exists an ordering $(\pi_i)_i$ of $P$ such that $d^{\pi_{i+1}}_t\leq \frac12 d^{\pi_i}_t$ for all $i$,
\item $|\bigsqcup_{\pi\in P}\cH^\pi|\geq \frac{|\cH|}{2q}$.
    \end{enumerate}
\end{definition}
Every uniform hypergraph on $[n]$ with sufficiently many hyperedges contains a weakly regular subpartition.
\begin{fact}\label{lem:pickt}
    Let $\cH$ be a $q$-uniform hypergraph on $[n]$ with $m\geq n^{2-2/q}$ total hyperedges. Then there exists a $t$-weakly regular subpartition of $\cH$ for some $t\in [q]$.
\end{fact}
\begin{proof}
    This follows from the decomposition algorithm of \cite[Lemma~4.7]{BHKL25} with $\eta\leq \frac12$. By Pigeonhole, there exists some $t\in [q]$ such that the union of all $t$-approximately strongly regular subhypergraphs contains at least $\frac{|\cH|}{2q}$ hyperedges, and the ordering property follows from the fact that if $d_t$ is the maximum degree of a $t$-tuple at a step where we construct a $t$-approximately strongly regular subhypergraph, we remove all $t$-tuples of weight $\geq d_t/2$.
\end{proof}
We now define weighted Kikuchi graphs, which are weighted unions of the Kikuchi graphs from \Cref{def:kikuchi} generated from each $t$-approximately strongly regular part of a weak regularity subpartition.

\begin{definition}[Weighted Level-$\ell$ Kikuchi Graph] \label{def:weightedkikuchi}
Suppose we have integers $1\leq t\leq q$ and $q$-uniform hypergraph matchings $\cH = (\cH_1,\dots,\cH_k)$ on $[n]$, and a subpartition $\cP=\{\cH^\pi:\pi\in P\}$ of disjoint $t$-approximately strongly regular subhypergraphs of $\cH$. For each $\pi\in P$, let $d^\pi_t$ be the $t$-degree of $\cH^\pi$ and let $K^\pi=K_\ell(\cH^\pi)$ be the level-$\ell$ Kikuchi graph of $\cH^\pi$. We define the level-$\ell$ Kikuchi graph $K_\ell(\cH,\cP)\coloneqq \bigcup_{\pi\in P}\frac{1}{d^\pi}K^\pi_\ell(\cH,\cP)$ to be the weighted union of the unweighted graphs $K^\pi=K^\pi_\ell(\cH,\cP)$, where every $e\in E(K^\pi)$ is given weight $w(e)=1/d^\pi$ in $K_\ell(\cH,\cP)$.

For all $\pi \in P$ and $i\in [k]$, let $\deg^\pi_i(S)\coloneqq \deg^{\cH^\pi}_i(S)$ and $\deg^\pi(S)\coloneqq \sum_{i\in [k]}\deg_i^\pi(S)$. We let $\deg_i(S)\coloneqq\sum_{\pi\in P}\frac{1}{d^\pi}\deg_i^\pi(S)$ and $\deg(S)\coloneqq\sum_{\pi\in P}\deg_i(S)$ to be the (weighted) degrees in $K_\ell(\cH,\cP)$.

\end{definition}

\newcommand{\degp}{\deg^{LR}}
\newcommand{\degtp}{\deg^{\tau,LR}}
\newcommand{\degpp}{\deg^{\pi,LR}}

When $\cH$ has sufficiently many hyperedges, the weighted Kikuchi graph $K_\ell(\cH,\cP)$ contains a rainbow cycle.
\begin{restatable}{lemma}{kikuchirainbowcycle}
    \label{lem:kikuchirainbowcycle}
    Suppose we have integers $1 \leq t \leq q$ where $q$ is odd and $q$-uniform hypergraph matchings
$\cH = (\cH_1, \cdots, \cH_k)$ such that $|\cH_i|\geq \delta n$ for all $i\in [k]$ and a $t$-weakly regular subpartition
$\mathcal{P} = \{\cH^\pi : \pi \in P\}$ of $\cH$. Suppose that
$k\geq q^{O(q)}\delta^{-1} n^{1-2/q} \log^2 n.$
Then the level-$\ell=n^{1-2/q}$ Kikuchi graph $K=K_\ell(\cH,\mathcal{P})$ contains a rainbow cycle.
\end{restatable}
Using this, we can now prove \Cref{lem:rainbowbound}.
\begin{proof}[Proof of \Cref{lem:rainbowbound}]
    By \Cref{lem:kikuchirainbowcycle}, there exists a cycle $V_1, \cdots, V_s, V_{s+1}=V_1$ in $K_\ell(\cH,\cP)$ such that $V_u\oplus V_{u+1}=C_u\oplus C'_{u}$ for some $C_u\in \cH_{i_u},C'_u\in \cH_{i'_u}$ and all $i_u,i'_u$ are distinct. We have
    \begin{align*}
        \left(\bigoplus_{u=1}^s C_u\right) \oplus \left(\bigoplus_{u=1}^s C'_u\right)=\bigoplus_{u=1}^s (V_u\oplus V_{u+1})=\emptyset\mcom
    \end{align*}
    so $\krain{q}\leq O(\delta^{-1}n^{1-2/q}\log^2 n)$.
\end{proof}
\begin{remark}
    While the Kikuchi graphs used here are based on those of \cite{BHKL25}, we note that the same result can be obtained using the Kikuchi graphs from \cite{JM25}.
\end{remark}
\subsection{Rainbow cycles in weighted Kikuchi graph: Proof of \Cref{lem:kikuchirainbowcycle}}
We first state the following property of Kikuchi graphs.
\begin{fact}[Observation 3.7 of \cite{HsiehKM23}] \label{fact:bound-D}
    The Kikuchi graph defined in \Cref{def:kikuchi} has $N = \binom{2n}{\ell}$ vertices, and each pair $(C,C')$ contributes a $D_\ell$-sized matching in the Kikuchi graph, where
    \begin{equation*}
        D_\ell \coloneqq \binom{q-t}{\floor{\frac{q-t}{2}}}^2  \binom{2n-2(q-t)}{\ell-(q-t)}  \cdot 2^{\1(\text{$q-t$ is odd})}
        \geq N \parens*{\frac{\ell}{2n}}^{q-t} \mper
    \end{equation*}
\end{fact}
Using this, we prove the following lower bound on the degree of a weighted Kikuchi graph.
\begin{lemma}\label{lem:avgdegbound}
    Let $\cP$ be a $t$-weakly regular subpartition of $\cH$ and $K_\ell=K_\ell(\cH,\cP)$ be the weighted Kikuchi graph defined in \Cref{def:weightedkikuchi}. Then the average degree $\bar{d}(K_\ell)$ of $K_\ell$ satisfies
    $\bar{d}(K_\ell)\geq \frac{|\cH|D}{48qN}$.
\end{lemma}
\begin{proof}
    For all $\pi\in P$, the number of parts $p_t^\pi$ in the $t$-approximately strongly regular decomposition $\cH^\pi=\bigsqcup_{\theta\in [p^\pi_t]}\chtp$ satisfies $\frac{|\cH^\pi|}{d^\pi_t}\leq p_t^\pi\leq \frac{2|\cH^\pi|}{d^\pi_t}$. Then
    \begin{align*}
        \bar{d}(K_\ell)=\sum_{\pi\in P}\frac{1}{d^\pi_t}\bar{d}(K^\pi)
        =\frac1N\sum_{\pi\in P}\frac{1}{d^\pi_t}\cdot D_\ell\sum_{\theta\in [p_t^\pi]}{|\chtp|\choose 2}
        &\geq \frac{D_\ell}{N}\sum_{\pi\in P}\frac{1}{d^\pi_t}\frac{|\cH^\pi|}{d^\pi_t}{d^\pi_t/2 \choose 2}
        \\&\geq \frac{D_\ell}{24N}\sum_{\pi\in P}|\cH^\pi|
        \\&\geq \frac{|\cH|D_\ell}{48 qN}\mper
    \end{align*}
\end{proof}
To prove \Cref{lem:kikuchirainbowcycle}, we will find a subgraph of the weighted Kikuchi graph satisfying the conditions of \Cref{thm:fullrainbow}. We define the left-right partitioned graph based on a partition of the $\cH_i$s, which allows us to ensure that given any edge $(S,T)$ induced by some $C\in\cH$, once we know which $i\in [k]$ such that $C\in \cH_i$, we can identify whether the elements of $C$ in $S\oplus T$ are $C\times \{1\}$ or $C\times \{2\}$.
\begin{definition}
    Let $K_\ell=K_\ell(\cH,\cP)$ be a weighted Kikuchi graph and $L\sqcup R=[k]$ be a partition of $[k]$. We define the weighted $L$-$R$ partitioned Kikuchi graph $K_{LR}$ to be the subgraph of $K_\ell$ such that an edge $(S,T)$ with $\gamma(S,T)=(i,j)$ is kept if $i\in L$ and $j\in R$. For all $\pi\in P$, we define $K^\pi_{LR}$ similarly. For all $\degpp_i(S)$ be the number of edges in $K^\pi_{LR}$ containing color $i$ incident to $S$ and $\degpp(S)\coloneqq \sum_{i\in [k]}\degpp(S)$, and $\degp_i(S)\coloneqq \sum_{\pi\in P}\frac{1}{d^\pi_t}\degpp_i(S)$ and $\degp(S)\coloneqq\sum_{\pi\in P}\degp_i(S)$.
\end{definition}
\begin{proposition}\label{prop:trivial}
    There exists a partition $L\sqcup R$ such that $\bar{d}(K_{LR})\geq\frac14 \bar{d}(K_\ell)$.
\end{proposition}
\begin{proof}
    Consider a uniformly random partition $L'\sqcup R'=[k]$. Then each edge is retained with probability $\frac14$, so there exists a partition $L\sqcup R=[k]$ such that $\bar{d}(K_{LR})\geq \frac14 \bar{d}(K_{\ell})$.
\end{proof}


\begin{proof}[Proof of \Cref{lem:kikuchirainbowcycle}]
\newcommand{\lqt}{\left\lfloor\frac{q-t}{2}\right\rfloor}
\newcommand{\cqt}{\left\lceil\frac{q-t}{2}\right\rceil}
    Let $A$ such that $k=A\delta^{-1}n^{1-2/q}\log^2 n$ and define $\Delta\coloneqq \frac{\bar{d}(K_\ell)}{32000\ell\log^2  n}$. We will show that there exists a subgraph $K'$ of $K_\ell(\cH,\cP)$ such that $d(K')\geq$ and for each $S\in V(K')$ and for any color $c$, the sum of the weights $w(e)$ of all edges $e$ incident to $S$ such that $\gamma(e)\ni c$ is at most $\Delta$.
By \Cref{lem:avgdegbound,fact:bound-D}, we have 
\begin{align*}
    \bar{d}(K_\ell)
    &\geq\frac{|\cH|D_\ell}{24qN}
    \geq \frac{k\delta n}{24q}\left(\frac{\ell}{2n}\right)^{q-t}
    \geq A n^{2-2/q}\log^2 n\cdot n^{-2+2t/q}
    =An^{-1+2t/q}\ell\log^2 n\mcom
\end{align*}
so $\Delta\geq An^{-1+2t/q}$.
By \Cref{prop:trivial}, there exists a partition $L\sqcup R=[k]$ such that $\bar{d}(K_{LR})\geq\frac14 \bar{d}(K_\ell)$. Delete every edge $(S,T)$ such that $i\in \gamma(S,T)$ and $\degp_i(S)>\Delta$, and let $K'$ be the resulting graph. 

We now bound the number of edges we delete. Fix some $\pi\in P$, and consider a uniformly random edge $(S,T)\in K^\pi_{LR}$, in other words, some $S\xleftrightarrow{C, C'} T$ where $C,C'\in \cH^\pi$ and $C\in \cH_i,C'\in\cH_j$ where $i\in L,j\in R$. Let $\gamma(S,T)=\{i,j\}$, so $C\in \cH_i$ and $C'\in \cH_j$. Then we can bound the total weight deleted by
\begin{align*}
        \sum_{\pi\in P}\frac{D}{d^\pi_t}&\sum_{C,C'\in \cH^\pi}\Pr_{\substack{S\xleftrightarrow{C, C'} T\\C,C'\in \cH^\pi}} \bracks*{(S,T) \text{ is removed}}
        \\&\leq\sum_{\pi\in P}\frac{D}{d^\pi_t}\sum_{C,C'\in \cH^\pi}\frac{1}{\Delta} \E_{\substack{S\xleftrightarrow{C, C'} T\\C,C'\in \cH^\pi}}\bracks*{\max\left(\degp_i(S),\degp_j(S),\degp_i(T),\degp_j(T)\right)}
        \\&\leq\sum_{\pi\in P}\frac{D}{d^\pi_t}\sum_{C,C'\in \cH^\pi}\frac{2}{\Delta} \E_{\substack{S\xleftrightarrow{C, C'} T\\C,C'\in \cH^\pi}}\bracks*{\degp_i(S)+\degp_j(S)}
        \\&\leq2D\sum_{\pi\in P}\sum_{C,C'\in \cH^\pi}\frac{1}{\Delta d^\pi_t} \E_{\substack{S\xleftrightarrow{C, C'} T\\C,C'\in \cH^\pi}}\bracks*{\sum_{\tau\in P}\frac{1}{d_t^\tau}\left(\degtp_i(S)+\degtp_j(S)\right)}    \\&\leq4D\sum_{\pi\in P}\sum_{C,C'\in \cH^\pi}\frac{1}{\Delta d^\pi_t} \E_{\substack{S\xleftrightarrow{C, C'} T\\C,C'\in \cH^\pi}}\bracks*{\sum_{\substack{\tau\in P\\ \tau\geq \pi}}\frac{1}{d_t^\tau}\left(\degtp_i(S)+\degtp_j(S)\right)}\mper
\end{align*}
We now bound this quantity by considering all edges $(S,T')=S\xleftrightarrow{C_1, C_2} T'$ (where $C_1,C_2\in \cH^\tau$) incident to $(S,T)=S\xleftrightarrow{C, C'} T$ (where $C,C'\in \cH^\pi$) at $S$ which contribute to $\degtp_i(S)$ for some $i$. We do casework on $|\{C,C'\}\cap\{C_1,C_2\}|$ and $|C'\cap C_2|$. Note that this quantity is either $0$ or $1$, since $\{C,C'\}\neq \{C_1,C_2\}$ as the $\cH_i$s are matchings. 

\begin{itemize}
    \item \textbf{Case 1: $|\{C,C'\}\cap\{C_1,C_2\}|=1$.} We show that for every $C,C'$, the total contribution to $\frac{1}{\Delta} \E_{\substack{S\xleftrightarrow{C, C'} T\\C,C'\in \cH^\pi}}\bracks*{\sum_{\tau\in P}\frac{1}{d_t^\tau}\degtp_i(S)}$ due to some edge $(S,T)=S\xleftrightarrow{C, C'} T$ being incident to an edge $(S,T')=S\xleftrightarrow{C_1, C_2} T'$ that shares one color with $(S,T)$ is small. Given $C,C'$, when choosing $C_1,C_2$ such that $|\{C,C'\}\cap\{C_1,C_2\}|=1$, without loss of generality (and losing a factor of $2$), assume that $C_1=C$ and $C_2=C''\notin \{C,C'\}$.  First, note that this implies $C''\in \cH^\pi$, and letting $\theta$ such that $C\in \chtp$, we see that $C',C''\in \chtp$, so all nonzero contributions to this case occur when $\tau=\pi$. Also, $Q_\theta\in C',C''$ and thus $|C'\cap C''|\geq t$. First, choose some $a\geq 0$ such that $|(C'\cap C''|=t+a$. Then note that there are at most ${q\choose t+a}d^\pi_{t+a}\leq q^qd^\pi_{t+a}$ choices of $C''\in \cH^\pi$ by $t$-approximate strong regularity. 

    Since $S\xleftrightarrow{C, C''} T'$, $|S\cap ((C''\setminus Q_\theta)\times \{2\})|\geq \lqt$, and thus $|S\cap ((C''\setminus C')\times \{2\})|\geq\lqt-a$. 
    However, $(S\setminus T)\subset  ((C\setminus Q_\theta)\times \{1\})\cup ((C'\setminus Q_\theta)\times \{2\})$, which implies that $S\cap T$ must include the at least $\lqt-a$ elements of $S\cap ((C''\setminus C')\times \{2\})$. Union bounding over the choice of $\lfloor\frac{q-t}{2}\rfloor-a$ elements and since $S\cap T$ is a uniformly random element of ${[n]\times [2]\choose \ell-(q-t)}$, the probability this occurs is 
    \begin{align*}
        {q-t-a\choose \max\left(\lqt-a,0\right)}\frac{{2n-\max\left(\lqt-a,0\right)\choose \ell-(q-t)-\max\left(\lfloor\frac{q-t}{2}\rfloor-a,0\right)}}{{2n\choose \ell-(q-t)}}\leq O(q^q)\cdot \left(\frac{\ell}{n}\right)^{\max\left(\lqt-a,0\right)}\mper
    \end{align*}
    Therefore, for any $C,C'$, the contribution of all edges incident to some $(S,T')=S\xleftrightarrow{C_1, C_2} T'$ incident to $S\xleftrightarrow{C, C'} T$ at $S$ to the expectation term for this case is
    \begin{align*}
        \frac{4}{\Delta}\frac{1}{d^\pi_t}\sum_{a=0}^{q-t}2d^\pi_{t+a}\cdot O(q^q)\left(\frac{\ell}{n}\right)^{\max\left(\lqt-a,0\right)}
        &\leq \frac{O(q^q)}{A}\left(\cqt+\sum_{a=0}^{\lqt}\frac{d_{t+a}^\pi}{d_t^\pi}n^{1-\frac{2t}{q}-\frac2q(\lqt-a)} \right)
        \\&\leq \frac{O(q^q)}{A}\left(1+\sum_{a=0}^{\lqt}\frac{d_{t+a}^\pi}{d_t^\pi}n^{\frac{2a}{q}-\frac1q(t-\1\{t\text{ even}\})}\right)
        \\&\leq \frac{O(q^q)}{A}\mcom
    \end{align*}
    where the last inequality follows from the second condition of \Cref{def:goodindex}. Summing over $\pi\in P$ and $C,C'\in \cH^\pi$, we delete an $O(q^q)/A$ fraction of the weight induced by incident edges such that $|\{C,C'\}\cap \{C_1,C_2\}|=1$.
    
    \item \textbf{Case 2: $|\{C,C'\}\cap\{C_1,C_2\}|=0$ and $C'\cap C_2=\emptyset$.}
Again, we show that for every $C,C'$, the total contribution to $\frac{1}{\Delta} \E_{\substack{S\xleftrightarrow{C, C'} T\\C,C'\in \cH^\pi}}\bracks*{\sum_{\tau\in P}\frac{1}{d_t^\tau}\degtp_i(S)}$ due to some edge $(S,T)=S\xleftrightarrow{C, C'} T$ being incident to an edge $(S,T')=S\xleftrightarrow{C_1, C_2} T'$ that shares one color with $(S,T)$ is small.
    Given $C,C'$, without loss of generality, assume that the shared color is $C,C_1\in \cH_i$ for some $i\in [k]$. Then there are at most $|\cH_i\cap \cH^\tau|$ choices to choose $C_1\in \cH_i$, which fixes $\sigma$ such that $C_1\in \cH^\tau_\sigma$ and $Q_\sigma\in C_1,C_2$, so there are now $d_t^\tau$ choices for $C_2$ in $\cH^\tau_\sigma$.

    Since $S\xleftrightarrow{C_1, C_2} T'$, $|(S\setminus T')\cap ((C_1\setminus Q_\sigma)\times \{1\})|\coloneqq s\in \{\lqt,\cqt\}$. Then since $C\cap C_1=\emptyset$, $S\cap T$ must include these $s$ elements. Furthermore, we have $|S\cap ((C_2\setminus C'\setminus Q_\sigma)\times \{2\})|=|(S\setminus T')\cap ((C_2\setminus Q_\sigma)\times \{2\})|= q-t-s$, and thus $S\cap T$ must contain these $q-t-s$ elements. Thus we can first choose these $s$ elements of $(C_1\setminus Q_\sigma)\times\{1\}$ and $q-t-s$ elements of $(C_2\setminus Q_\sigma)\times\{2\}$, and the (uniformly distributed) $S\cap T$ must contain these $q-t$ total elements. Union bounding over the choice of $s$ and $q-t-s$ total elements from $C_1\setminus Q_\sigma$ and $C_2\setminus Q_\sigma$, the probability this occurs is at most
    \begin{align*}
        {q-t\choose s}{q-t\choose q-t-s}\frac{{2n-(q-t) \choose\ell-2(q-t)}}{{2n\choose\ell-(q-t)}}\leq O(q^{q})\left(\frac{\ell}{n}\right)^{q-t}\mper
    \end{align*}
    Thus the contribution of all edges incident to some $(S,T')=S\xleftrightarrow{C_1, C_2} T'$ incident to $S\xleftrightarrow{C, C'} T$ at $S$ to the expectation term for this case is
    \begin{align*}
        \frac{4}{\Delta }\sum_{\substack{\tau\in P\\\tau\geq \pi}}\frac{1}{d^\tau_t}|\cH_i\cap\cH^\tau|d_t^\tau\cdot \left(\frac{\ell}{n}\right)^{q-t}
        &=\frac{1}{A}\sum_{\substack{\tau\in P\\\tau\geq \pi}}n^{1-\frac{2t}{q}}|\cH_i\cap\cH^\tau|n^{-\frac{2}{q}(q-t)}
        \\&\leq \frac{1}{A}\sum_{\substack{\tau\in P\\\tau\geq \pi}}\frac{|\cH_i\cap \cH^\tau|}{n}
        \leq \frac{|\cH_i|}{An}\leq \frac{1}{A}\mper
    \end{align*}
    
    \item \textbf{Case 3: $|\{C,C'\}\cap\{C_1,C_2\}|=0$ and $C'\cap C_2\neq\emptyset$.} For this case, we bound the entire sum rather than the expectation term. We count the number of ways to pick $C,C'\in \cH^\pi$ and $C_1,C_2\in \cH^\tau$. First, choose some $a> 0$ such that $|C'\cap C_2|=a$. Without loss of generality, suppose that the shared color comes from $C,C'$, so there exists some $i\in [k]$ such that $C,C_1\in \cH_i$ and we have $C\cap C_1=\emptyset$. Given $\pi\in P$, there are $|\cH^\pi|$ ways to pick $C'\in \cH^\pi$, which determines $\theta$ such that $C,C'\in \cH^\pi_\theta$ and $Q_\theta\subset C,C'$. Then there are ${q\choose a}d^\tau_a\leq q^ad^\tau_a$ ways to choose $C'\cap C_2$ and then $C_2$. This fixes $\sigma$ such that $C_1,C_2\in \cH^\tau_\sigma$ and $Q_\sigma\in C_1,C_2$.
    It now remains to pick $C,C_1$. We choose $C_1$ first, and there are at most $d_t^\tau$ possible such choices such that $ Q_\sigma\subset C_1$. This determines $i\in [k]$ such that $C_1\in \cH_i$, and then by the matching property there is at most one $C\in \cH_i$ such that $Q_\theta\subset C$. Thus there are $q^a|\cH^\pi|d^\tau_a d_t^\tau$ ways to pick $C,C',C_1,C_2$.

        Since $S\xleftrightarrow{C_1, C_2} T'$, $|(S\setminus T')\cap ((C_1\setminus Q_\sigma)\times \{1\})|\coloneqq s\in \{\lqt,\cqt\}$. Then since $C\cap C_1=\emptyset$, $S\cap T$ must include these $s$ elements. Furthermore, we have $|(S\setminus T')\cap ((C_2\setminus Q_\sigma)\times \{2\})|= q-t-s$, so $|S\cap ((C_2\setminus C'\setminus Q_\sigma)\times \{2\})|\geq q-t-s-a$, and thus $S\cap T$ must contain these $q-t-s-a$ elements. Thus we can first choose these $s$ elements of $(C_1\setminus Q_\sigma)\times\{1\}$ and $q-t-s-a$ elements of $(C_2\setminus Q_\sigma)\times\{2\}$, and the (uniformly distributed) $S\cap T$ must contain these $s+\max(q-t-s-a,0)=\max(q-t-a,s)\geq \max(q-t-a,\lqt)$ elements. The probability this occurs is at most
    \begin{align*}
        {q-t\choose \lfloor\frac{q-t}{2}\rfloor}{q-t\choose \max\left(\lceil\frac{q-t}{2}\rceil-a,0\right)}\frac{{2n-\max(q-t-a,\lfloor\frac{q-t}{2}\rfloor) \choose\ell-(q-t)-\max(q-t-a,\lfloor\frac{q-t}{2}\rfloor)}}{{2n\choose\ell-(q-t)}}\leq O(q^{2q})\left(\frac{\ell}{n}\right)^{\max(q-t-a,\lfloor\frac{q-t}{2}\rfloor)}\mper
    \end{align*}
     Therefore, we can bound the contribution to the deleted weight by
    \begin{align*}
        4D &\sum_{\pi\in P}\frac{2}{\Delta d^\pi_t}\sum_{\substack{\tau\in P\\ \tau\geq \pi}}\frac{1}{d^\tau_t}O(q^{2q})\sum_{a=1}^q |\cH^\pi|d_a^\tau d_t^\tau\left(\frac{\ell}{n}\right)^{\max(q-t-a,\lfloor\frac{q-t}{2}\rfloor)}
        \\&\leq \frac{D\cdot  O(q^{2q})}{\Delta}\sum_{\pi\in P}|\cH^\pi|\sum_{\substack{\tau\in P\\ \tau\geq \pi}}\sum_{a=1}^q \frac{d_a^\tau}{d_t^\pi}\left(\frac{\ell}{n}\right)^{\max(q-t-a,\lfloor\frac{q-t}{2}\rfloor)}
        \\&\leq \frac{D\cdot O(q^{2q})}{A}\sum_{\pi\in P}|\cH^\pi|\sum_{\substack{\tau\in P\\ \tau\geq \pi}}n^{1-\frac{2t}{q}}\left(\cqt \frac{d_{\cqt}^\tau}{d^\pi_t}n^{-\frac1q(q-t-\1\{t\text{ even}\})}+\sum_{a=1}^{\lqt} \frac{d_a^\tau}{d_t^\pi}n^{-\frac{2}{q}(q-t-a)}\right)
        \\&\leq \frac{D\cdot O(q^{2q})}{A}\sum_{\pi\in P}|\cH^\pi|\sum_{\substack{\tau\in P\\ \tau\geq \pi}}\left( \frac{d_{\cqt}^\tau}{d^\pi_t}n^{-\frac{t-\1\{t\text{ even}\}}{q}}+\frac{d^\tau_t}{d^\pi_t}\sum_{a=1}^{\lqt}\frac{d^\tau_a}{d^\tau_t}n^{-1+\frac{2a}{q}}\right)
        \\&\leq \frac{D\cdot O(q^{2q})}{A}\sum_{\pi\in P}|\cH^\pi|\left(1+\sum_{\substack{\tau\in P\\ \tau\geq \pi}}\frac{d^\tau_t}{d^\pi_t}\right)
        \\&\leq \frac{D\cdot O(q^{2q})}{A}\sum_{\pi\in P}|\cH^\pi|\leq D|\cH|\cdot \frac{O(q^{2q})}{A}\mcom
    \end{align*}
    where the fourth inequality follows from the first condition of \Cref{def:goodindex} and the fifth follows from the second condition of \Cref{def:weakregdef}. Thus we delete an $O(q^{2q})/A$ fraction of the weight induced by incident edges in this case. 
    
    \end{itemize}
    Thus, by picking $A$ to be sufficiently large, in each of the three cases, we delete at most a $\frac{1}{24}$ fraction of edges in $K_\ell$, so $\bar{d}(K')\geq\frac12\bar{d}(K_{LR})\geq\frac18\bar{d}(K_\ell)\geq 16000\Delta \log N\log\log N$. By \Cref{thm:fullrainbow}, $K'$ (and therefore $K_\ell$) contains a rainbow cycle.
    \end{proof}

\section{Acknowledgements}
We thank Omar Alrabiah and Pravesh K. Kothari for helpful discussions, and to Arpon Basu, Matija Buci\'c, and Jun-Ting (Tim) Hsieh for helpful discussions on an earlier related project.
\section{AI Statement}
A weaker version of \Cref{thm:fullrainbow} for unweighted properly list-edge colored graphs with list-size $2$ was first obtained with assistance from ChatGPT 5.6-Sol, after which \Cref{thm:lccbound} was obtained by the author with minimal AI use. \Cref{thm:fullrainbow} was later strengthened to its current version using ChatGPT 6-Astra. In addition, ChatGPT 6-Astra was used for literature review and proofreading. The paper is entirely manually written and the author takes full responsibility for the correctness and attribution of all content.
\bibliographystyle{alpha}
\bibliography{bib}

\end{document}